\documentstyle[12pt]{article}
\def\bea{\begin{eqnarray}}
\def\beq{\begin{equation}}
\def\eea{\end{eqnarray}}
\def\eeq{\end{equation}}

\newcommand{\app}[1]{\setcounter{section}{0}
\setcounter{equation}{0}\renewcommand{\thesection}
{\Alph{section}}\section{#1}}

\begin{document}
\begin{center}
{\Large\bf Some examples of classical coboundary Lie bialgebras
with coboundary duals}
\end{center}
\begin{center}
\bf M A Sokolov
\end{center}
\begin{center}
St Petersburg Institute of Machine Building,
Poliustrovskii pr 14, 195108, St Petersburg, Russia\\
E-mail addresses: sokolov@pmash.spb.su \, \& \, mas@ms3450.spb.edu
\end{center}

\begin{center}
Abstract
\end{center}
Some examples are given of finite dimensional Lie bialgebras whose
brackets and cobrackets are determined by pairs of $r$-matrices.
\vspace{1cm}

The aim of this Letter is to give some low-dimensional
examples of classical coboundary Lie bialgebras \cite{D,CP}
with coboundary duals. Since such structures can be specified (up to
automorphisms) by pairs of $r$-matrices, so it is natural to
call them {\it bi-$r$-matrix bialgebras} (B$r$B).
There are some reasons to study B$r$B.
   Assuming that both Lie bialgebras of a dual pair are coboundary, we
impose additional constrains, which can facilitate the search of
new classical $r$-matrices connected with nonsemisimple Lie
algebras. Recall that many Lie algebras of physical interest are
nonsemisimple ones, but up to now there is detailed classification
of $r$-matrices only for the complex simple Lie algebras \cite{BD}.
The presence of a pair of $r$-matrices is useful for practical
quantization of Lie bialgebras \cite{D} permitting more symmetrical
treatment of quantum algebras and groups.
     However, most interesting applications of B$r$B are possible in
the theory of bihamiltonian dynamical systems \cite{M,GD}. In this
case the presence of a pair of $r$-matrices allows us to define
the pair of dynamical systems on the space which is the space of
the original Lie algebra canonically identified with its dual space
\cite{S}.

Now recall some basic definitions \cite{D,CP}.
Let $L$ be a finite-dimensional Lie algebra and $L^{*}$ the dual
space of $L$ in respect to a nondegenerate bilinear form $<.,.>$
on $L^{*}\times L$.
The Lie algebra $L$ is called a bialgebra, if there exist a map
$\delta :L\rightarrow L\otimes L$ which is an 1-cocycle

\beq
\label{in3}
\delta ([x,y])=[\delta (x),y\otimes 1+1\otimes y]+[x\otimes
1+1\otimes x,\delta (y)], \quad x,y\in L
\eeq
and which defines on $L^{*}$ a Lie algebra structure
$[.,.]_{*}:L^{*}\otimes L^{*}\rightarrow L^{*}$ by the following
relation

\beq
\label{in4}
<[\xi ,\eta ]_{*},x>=<\xi \otimes \eta ,\delta (x)>,\quad
x\in L,\quad \xi ,\eta \in L^{*}.
\eeq
If one puts

\beq
\label{in5}
\delta (x)=[x\otimes 1+1\otimes x,r], \quad
r=r^{mn}e_m\otimes e_n\in L\otimes L
\eeq
then the 1-cocycle condition (\ref{in3}) is fulfilled identically.
In this case the Lie bialgebra $L$ is called a coboundary one.
It can be shown \cite{CP} that the Jacobi identity for the elements
of $L^{*}$ is equivalent to the following equation

\beq
\label{in8}
[x\otimes 1\otimes 1 + 1\otimes x\otimes 1 +
1\otimes 1\otimes x,[r,r]_s]=0, \quad x\in L.
\eeq
The Schouten brackets $[r,r]_s$ in the above formula are defined by

$$
[r,r]_s\equiv [r_{12},r_{13}+r_{23}]+[r_{13},r_{23}]]
$$
where, as usual, $r_{12}=r^{mn}e_m\otimes e_n\otimes 1,$ etc.
The simplest way to satisfy the relation (\ref{in8}) is to put

$$
[r,r]_s=0.
$$
This is the classical Yang-Baxter equation and its
solutions for the complex simple Lie algebras are classified in
the mentioned work \cite{BD}. If some $r$-matrix satisfies CYBE, then
the related bialgebra (as well as this $r$-matrix itself) is called a
quasitriangular bialgebra. If, in addition, this $r$-matrix
satisfies the unitarity condition

$$
\sigma \circ r = - r
$$
where $\sigma$ is the flip operator
$\sigma (x\otimes y)=y\otimes x$,
then the bialgebra is called a triangular one. When we leaving out
${\bf Z}_2-$graded Lie bialgebras, we can restrict ourselves by
unitary $r$-matrices because brackets $[.,.]_{*}$ are defined
only by the antisymmetric part of $r.$

Now extend slightly the above definition of bi-$r$-matrix
bialgebras. Let $L$ be a coboundary Lie bialgebra with an 1-cocycle
(\ref{in5}) and $L^{*}$ be its dual space endowed with a structure
of coboundary Lie bialgebra as well. Let $\delta ^{*}$ be a
corresponding 1-cocycle defined by some $r$-matrix $r^{*}$

$$
\delta ^{*}(\xi )=[\xi \otimes 1+1\otimes \xi ,r^{*}]_{*},
\quad \xi \in L^{*},\quad r^{*}\in L^{*}\otimes L^{*}.
$$
We will call the pair $(L,L^{*})$ a bi-$r$-matrix
bialgebra if the Lie brackets $[.,.]^{^{\prime }}$ on $L$
defined by $\delta ^{*}$

$$
<\delta ^{*}(\xi ),x\otimes y>=<\xi ,[x,y]^{^{\prime }}>,
\qquad x,y\in L,\qquad \xi \in L^{*}
$$
are equivalent to the original ones

$$
[x,y]^{^{\prime }}=S^{-1}[Sx,Sy],\quad x,y\in L,\quad S\in
{\rm Aut}(L).
$$
Because of the equal status of the algebras $L$ and $L^{*}$ in our
definition, we will call B$r$B under consideration a double
name in what follows.

The simplest and well known example of a self-dual B$r$B gives the case
of the two dimensional non-Abelian Lie algebra with the generators
$e_1,e_2$ \cite{D}. Let the normalized commutation relation
of this algebra be

\beq
\label{in12}
[e_1,e_2]=e_1.
\eeq
The only possible triangular $r$-matrix

$$
r=e_1\wedge e_2
$$
produces by (\ref{in4}), (\ref{in5}) the relation among
the generators $e^1, e^2$ of $L^{*}$ (we use the form
$<e^i,e_j>={\delta}^i_j$)

$$
[e^1,e^2]_{*}=e^2.
$$
The $r$-matrix
$r^{*}\in L^{*}\wedge L^{*}$

$$
r^{*}=-e^1\wedge e^2
$$
induces on $L$ the original brackets (\ref{in12}).

Consider the case of a three-dimensional complex Lie algebra $L$.
It is well known that the commutation relations among its generators
$\left\{e_i\right\},\,i=1,2,3$ can be reduced by a linear
transformation to the following form

\beq
\label{a1}
[e_1,e_2]=ae_2+be_3,\quad [e_2,e_3]=ce_1,\quad
[e_1,e_3]=fe_2+ae_3.
\end{equation}
The Jacobi identity restricts the values of the structure constants
$a,b,c,f$
by the condition $ac=0.$ At first consider the case
$a=0,c\neq 0.$ An arbitrary unitary $r$-matrix

\beq
\label{a4}
r=r^{12}e_1\wedge e_2+r^{13}e_1\wedge e_3+r^{23}e_2\wedge e_3
\eeq
gives the relations among the generators
$\left\{ e^i\right\},\,i=1,2,3$ of the dual Lie algebra $L^{*}$

\beq
\label{a5}
[e^1,e^2]_{*}=fr^{13}e^1-cr^{23}e^2,\quad
[e^1,e^3]_{*}=br^{12}e^1-cr^{23}e^3,\quad
[e^2,e^3]_{*}=br^{12}e^2-fr^{13}e^3
\eeq
which satisfy the Jacobi identity. Note that the structure constants and
$r$-matrix elements can be mutually absorbed
($fr^{13}\rightarrow r^{13}$ or $fr^{13}\rightarrow f,$ etc.).
For the Schouten brackets one has

\beq
\label{a6}
[r,r]_s=\left( b(r^{12})^2-f(r^{13})^2+c(r^{23})^2\right)
e_1\wedge e_2\wedge e_3.
\eeq
It is not difficult to show that any choice of the structure constants
$b,c,f$ and the $r$-matrix elements $r^{12},r^{13},r^{23}$ in the
non-Abelian case leads to commutators on $L^{*}$ which are equivalent
to the following ones

\beq
\label{a8}
[\tilde e^i,\tilde e^j]_{*}=\tilde e^j,\qquad
[\tilde e^j,\tilde e^k]_{*}=0,\qquad
[\tilde e^i,\tilde e^k]_{*}=\tilde e^k
\eeq
where $(i,j,k)$ is some transposition of (1,2,3). The Lie algebra
with the above commutation relations (we will denominate it $j(3)$)
appeared, in particular, in the semiclassical limit of the matrix quantum group $SL_q(2)$
and the Jordanian group $SL_h(2)$ (see, for example, \cite {VZW,MS}).

Taking into account the above results, consider the special case
$a=1,$ $b=c=f=0.$ For these structure constants the relations
(\ref{a1}) take the form

\beq
\label{a14}
[e_1,e_2]=e_2,\quad [e_2,e_3]=0,\quad [e_1,e_3]=e_3
\eeq
and the Schouten brackets for an arbitrary $r$-matrix (\ref{a4})
identically vanish $[r,r]_s=0.$ The commutation relations obtained
form (\ref{a4}),(\ref{a14}) by (\ref{in4}),(\ref{in5})

\beq
\label{a15}
[e^1,e^2]_{*}=r^{12}e^1,\quad [e^1,e^3]_{*}=r^{13}e^1,\quad
[e^2,e^3]_{*}=2r^{23}e^1-r^{13}e^2+r^{12}e^3.
\eeq
satisfy the Jacobi identity.

Summarizing the above formulae (\ref{a1}) -- (\ref{a15}),
list in the Table 1 (placed in the Appendix) the most interesting from
the physical point of view three-dimensional bi-$r$-matrix bialgebras:
the Heisenberg-Weyl algebra $h(3)$ ($a=b=f=0;c=1$) --- $j(3)$,
the Euclidean $e(2)$ ($a=b=0,f=-1;c=1$) --- $j(3)$ and the Poincar\`e
$p(2)$ ($a=b=0,f=-1;c=-1$) --- $j(3)$. For the sake of clearness we
consider the Lie algebras of the groups of the Euclidean and
pseudoeuclidean planes separately in spite of these cases are
distinguished unessentially in the structure constants and
coefficients of $r$-matrices. Remark, that some $r-$matrices defining
the brackets $[.,.]_*$ and presented in the Table 1 are well known.
For example, $r=-e_2\wedge e_3$ was considered in the papers
\cite{CGST,K,BHP,Ko}, and $r=e_1\wedge e_3+ie_2\wedge e_3$ in
\cite{So}.

Both of the dual algebras for $e(2)$ and $p(2)$, as was mentioned
above, are equivalent to the one with the commutation relations
(\ref{a8}). It can be displayed explicitly. For example, in the
case $p(2),$ using the $r$-matrix $r=e_2\wedge e_3$ instead of
$r=e_1\wedge e_3+e_2\wedge e_3,$ one obtains the
relations

$$
[e^1,e^2]_*=e^2,\quad [e^2,e^3]_*=0,\quad [e^1,e^3]_*=e^3
$$
directly. By the $r$-matrix $r^{*}=e^1\wedge e^3$ one has

$$
[e_1,e_2]=0,\quad [e_1,e_3]=e_1,\quad [e_2,e_3]=-e_2.
$$
After the transformation

$$
S:e_1\longmapsto \tilde e_1=e_1+e_2,\quad
e_2\longmapsto \tilde e_2=e_1-e_2,
e_3\longmapsto \tilde e_3=-e_3,
$$
we restore the original $p(2)$ commutation relations.
Let us stress here once more that all of the considered
three-dimensional bialgebras $h(3),e(2),p(2)$ and
the bialgebra $sl(2,C)$ have the same Lie algebra as a dual
counterpart. This fact can be used for their unified multiparameter
quantization \cite{DKS}.

Until now we used the unitary $r$-matrices because only
their antisymmetric parts $r^a=1/2(r-\sigma \circ r)$ define
the brackets for the dual Lie algebras. When we turn to the
case of Lie superalgebras the symmetric parts
$r^s=1/2(r+\sigma \circ r)$ play the important role.
We placed in the Table 1 the three-dimensional bi-$r$-matrix
superbialgebra in which the dual pair is
the Poincar\`e algebra $p(2)$ and the anticommutator algebra
$c(3)$. In this case the equivalence of the brackets $[.,.]^\prime$
given by the $r$-matrix $r^*=e^1\wedge e^2$ and the original $p(2)$ --
brackets is displayed by the transformation
$ S:\,
e^3\longmapsto \tilde e^3,\quad -e^1+e^2\longmapsto \tilde e^1,
\quad e^1+e^2\longmapsto \tilde e^2.
$

In fact, in the Table 1 are listed all possible three-dimensional B$r$B
and no one of them is self-dual. However, the extending of
the Heisenberg-Weyl algebra $h(3)$ gives the self-dual example. This is
so-called (1+1) extended  Galilei algebra. This algebra was studied
in detailed in the work \cite{BCH} and the formulae placed in the Table 1
are taken from there.

It is not difficult to generalize the above low-dimensional
examples to the $n-$dimensional cases. To facilitate the search of
new B$r$B it is useful to take into account the following
assertions. Let $L$ be a $n$-dimensional complex Lie algebra, and
$I$ its subalgebra. As a space $L$ can be represented by the sum
$L=I\oplus K.$ Denote by $I^{\prime }$ and $K^{\prime }$ the
subspaces of the dual space $L^{*} = I^{\prime}\oplus K^{\prime}$
which are orthogonal to $K$ and $I$ respectively:
$<I^{\prime},K>=<K^{\prime },I>=0.$ Let us list
simple properties which are the direct consequences of the formula
(\ref {in5}) provided that
$\delta (x)=[x\otimes 1+1\otimes x,r],$ $x\in L,$
$r\in L\otimes L$ defines on $L$ the structure of a coboundary Lie
bialgebra.

\begin{itemize}
\item  If $I$ is a subalgebra of $L,$ $r\in L\otimes L \ominus
K\otimes K ,$ then $K^{\prime}$ is an subalgebra of $L^{*}.$

\item  If $I$ is an ideal, then $K^{\prime }$ is
       a subalgebra of $L^{*}.$

\item  If $I$ is an ideal, $r\in I\otimes I,$ then
       $K^{\prime }$ is a centre of $L^{*}.$

\item  If $I$ is a centre, then $K^{\prime }$ is the commutant of
       $L^{*}.$

\item  If $I$ is a commutative ideal, $r\in I\otimes I,$
       then $K^{\prime}$ is both the commutant and the centre
       of $L^{*}.$

\item  If $I$ is both a commutant and a centre, $r\in K\otimes K,$
       then $K^{\prime}$ the commutative commutant.
\end{itemize}

The last two points allow directly generalize the Heisenberg-Weyl
bi-$r$-matrix bialgebra $h(3) - j(3)$ to the
$n-$dimensional case (see Table 1).

The author would like to thank E V Damaskinsky
and P P Kulish for discussions. This work was supported in
part by the Russian Fund for Fundamental Research (grants no
98-01-00310 and 00-01-00500).

\newpage
\renewcommand{\theequation}{\thesection.\arabic{equation}}
\app{Appendix}

\bigskip
{
\noindent
{\bf {Table 1}}.  Some bi-$r$-matrix bialgebras.

\medskip
\noindent
\begin{tabular}{|c|c|c|c|}
\hline
         $L$    & $r$  &   $L^*$    &      $r^*$ \\
\hline
\multicolumn{4}{|c|}{The Heisenberg-Weyl algebra $h(3) - j(3)$}\\
\hline

         $[e_1,e_2]=0  $ & &$[e^1,e^2]_*=e^2$ & \\

$[e_2,e_3]=e_1$ &$-e_2\wedge e_3$ \cite{CGST,K,BHP,Ko} &$[e^2,e^3]_*=0$&

$\frac 12e^2\wedge e^3$           \\

 $[e_1,e_3]=0$ &&$[e^1,e^3]_*=e^3$ &     \\
\hline
\multicolumn{4}{|c|}{The algebra of the Euclidean plane $e(2)
- j(3)$  }\\
\hline
$[e_1,e_2]=0 $ & &$[e^1,e^2]_*=-e^1-ie^2$ &  \\

$[e_2,e_3]=e_1$ &$e_1\wedge e_3+ie_2\wedge e_3$ \cite{So} &
$[e^2,e^3]_*=e^3$&
$-\frac 12(e^1\wedge e^3-ie^2\wedge e^3)$           \\

$[e_3,e_1]=e_2$ &&$[e^1,e^3]=-ie^3$ &\\
\hline
\multicolumn{4}{|c|}{The algebra of the pseudoeuclidean plane
$p(2) - j(3)$}\\
\hline

$[e_1,e_2]=0 $ & &$[e^1,e^2]_*=-e^1+e^2$ &  \\

$[e_2,e_3]=-e_1$ &$e_1\wedge e_3+e_2\wedge e_3$  &
$[e^2,e^3]_*=e^3$&
$-\frac 12(e^1\wedge e^3+e^2\wedge e^3)$           \\

 $[e_3,e_1]=e_2$ &&$[e^1,e^3]=e^3$ &
\\
\hline
\multicolumn{4}{|c|}{The algebra of the pseudoeuclidean plane $p(2)$ ---
The Clifford algebra $c(3)$}\\
\hline

$ [e_1,e_2]=0 $ & &$\{e^1,e^2\}_*=e^3$ & \\

$[e_3,e_2]=e_1$ &$1/2(e_1\otimes e_1+e_2\otimes e_2)$
& $[e^1,e^3]_*=0$ &  $e^1\wedge e^2$           \\

$[e_3,e_1]=e_2$ & & $[e^2,e^3]_* = 0$ &     \\
\hline

\multicolumn{4}{|c|}{Self-dual extended $(1+1)$ Galilei algebra}\\
\hline
$[e_3,e_1]=e_2$    &  & $[e^2,e^4]_*=e^1$  &           \\

$[e_3,e_2]=e_4$ & $e_1\wedge e_2 - e_3\wedge e_4$ \cite{BCH} &
$[e^1,e^4]_*=e^3$ &$-(e^1\wedge e^2 -e^3\wedge e^4) $\cite{BCH} \\

$[e_1,e_2]=0$ &  &$[e^1,e^2]_*=0$&
           \\

$[e_4,.]=0$ &&$[e^3,.]_*=0$ &     \\
\hline
\multicolumn{4}{|c|}{The Heisenberg-Weyl algebra $h(3n) - j(3n)$}\\
\hline
$[e_{2i},e_{3j}]=\delta _{ij}e_1$ &
                                  &
$[e^{2i},e^1]_{*}=a^ie^{2i}$  &        \\
$[e_{2i},e_{2j}]=0$ &
                                 &
$[e^{3i},e^1]_{*}=a^ie^{3i}$  &        \\

$[e_{3i},e_{3j}]=0$ &
$\sum_{i=1}^na^ie_{2i}\wedge e_{3i}$ \cite{CGST} &
$[e^{2i},e^{3j}]_{*}=0$ &
$-\frac 12\sum_{i=1}^n(a^i)^{-1}e^{2i}\wedge e^{3i}$ \\

$i,j=1,...,n $  &  & $i,j=1,...,n $ &\\
$[e_1,.]=0$& && \\
\hline
\end{tabular}

}

\end{document}